\documentclass[smallextended]{svjour3}
\usepackage{amsfonts}
\usepackage{latexsym}
\usepackage{amsmath}
\usepackage{amssymb}
\usepackage{amscd}
\usepackage{epsfig}
\usepackage[affil-it]{authblk}
\addtocounter{MaxMatrixCols}{4}
\newcommand{\undertilde}[1]{\ensuremath{\mathord{\vtop{\ialign{##\crcr
   $\hfil\displaystyle{#1}\hfil$\crcr\noalign{\kern1.5pt\nointerlineskip}
   $\hfil\tilde{}\hfil$\crcr\noalign{\kern1.5pt}}}}}}
\begin{document}

\title{Change Ringing and Hamiltonian Cycles : The Search for Erin and Stedman Triples}
\author{M. Haythorpe \and A. Johnson}
\institute{Michael Haythorpe
\at Flinders University\\
\email{michael.haythorpe@flinders.edu.au}
\and
Andrew Johnson
\at IBM Hursley\\
\email{andrew\_johnson@uk.ibm.com}
}

\maketitle {\abstract A very old problem in campanology is the search for peals. The latter can be thought of as a heavily constrained sequence of all possible permutations of a given size, where the exact nature of the constraints depends on which method of ringing is desired. In particular, we consider the methods of bobs-only Stedman Triples and Erin Triples; the existence of the latter is still an open problem. We show that this problem can be viewed as a similarly constrained (but not previously considered) form of the Hamiltonian cycle problem (HCP). Through the use of special subgraphs, we convert this to a standard instance of HCP. The original problem can be partitioned into smaller instances, and so we use this technique to produce smaller instances of HCP as well. We note that the instances known to have solutions provide exceptionally difficult instances of HCP.}

\section{Introduction and Terminology}\label{sec-Introduction}

Change ringing is the art of striking a set of bells in constantly varying orders, such that at each stage, each bell is struck once before the order changes. Typically, each bell has a designated person, or {\em ringer}, assigned to strike it, and at each stage, the order in which the ringers strike the bells depends on various schemes. Such an order is called a {\em row}, and initially, the bells are struck in decreasing order of the tone they produce; this row is called {\em rounds}. Each row can be represented as a sequence of $n$ numbers (assuming there are $n$ bells) where the first number denotes which bell should be struck first, the second number denotes which bell should be struck second, and so on. The bells are labelled in order from the highest tone to the lowest. For example, suppose there are five bells; then rounds could be represented as 12345. In this notation, we say that the first bell is in the first {\em place}, and so on.

After each row is struck, a {\em change} occurs to transition to a new row. However, changes are constrained such that each entry in the new row can only have moved, at most, one place from the previous row\footnote{This constraint is due to the physical limitations of ringing church bells. Rather than individual strikes, the bells are set in motion by manipulating long ropes, and the order of ringing is controlled by increasing or decreasing the speeds at which the various ropes are pulled. The weights of the bells make it difficult to alter the speed too significantly from one strike to the next, hence the constraint.}. It is not difficult to see that this means each bell must either retain its previous place, or swap places with a neighbour. For example, a valid change from rounds would be to transition to the row 13254, but no valid change transitions from rounds to 23451 since this would require the first bell to move four places. In cases where there is an odd number of bells, it is clear at least one bell must retain its place. A compact way of representing a sequence of changes is to simply list the places where the bell is unchanged (or list an 'x' if all bells change place), with each change separated by a dot; this is called {\em place notation}. The constraints on changes mean that the changes represented by any valid place notation are unique. For example, from rounds with five bells, the place notation 3.1.345.1 goes from rows 12345 to 21354 to 23145 to 32145 to 31254.

The rule governing the sequence of changes made is called the {\em method}, which is typically a short pattern which is learned by heart by the ringers. Ideally, the aim is to ring all possible rows before returning to rounds. The set of all possible rows is called the {\em extent}, which for $n$ bells will consist of $n!$ permutations of $n$ numbers. However, even for fairly small numbers of bells, there are no short patterns which can be repeated to achieve this. Hence, it is necessary for a conductor (who must also be one of the ringers) to make {\em calls}, which are alterations to the method, as necessary. Traditionally, as well as in modern competitive change ringing, neither the conductor nor the other ringers are permitted to bring in notes, so these calls must be memorised. Hence, it is often desirable if the calls occur at regular intervals, follow simple patterns, or are as infrequent as possible.

In this manuscript, we will focus on the case where there are seven bells; this case is called {\em triples}, so named because at most three swaps can occur per change. Note that although this implies there are only seven bells, in reality eight bells are used (corresponding to an octave), with the eighth and lowest bell always struck last in each row. Therefore, an extent for triples consists of $7! = 5040$ rows. For any number of bells $n \leq 7$, striking the full extent including the final return to rounds is called a {\em peal}; if $n \geq 8$ a peal must consist of at least 5000 changes. This number of changes will typically take several hours to complete.

For triples, there are only four ways for a change to induce three swaps. Specifically, any one of the first, third, fifth or seventh bells can stay in place. There are three common methods used for triples called Grandsire Triples, Stedman Triples and Erin Triples defined as follows. Grandsire Triples is the method induced by repeating the pattern 3.1.7.1.7.1.7.1.7.1.7.1.7.1. Stedman Triples is the method induced by repeating the pattern 3.1.3.1.3.7.1.3.1.3.1.7. Erin Triples is the method induced by repeating the pattern 3.1.3.1.3.7\footnote{In reality, ringers of Stedman Triples start in the tenth place, that is, they repeat the pattern 3.1.7.3.1.3.1.3.7.1.3.1, while ringers of Erin Triples start in the sixth place, repeating the pattern 7.3.1.3.1.3. However, given that a peal forms a cycle, it is mathematically equivalent to start at any place in the method, and it will be convenient in what follows to orient Stedman and Erin Triples so that both methods start off identically.}. None of these methods produce an extent; indeed, Grandsire Triples returns to rounds after 70 changes, Stedman Triples after 84 changes, and Erin Triples after 42 changes. However, it is possible to make calls to alter the methods appropriately, with two standard types of calls being used named {\em bobs} and {\em singles}.

In the case of Grandsire Triples, the pattern consists of fourteen changes, where the thirteenth change is to keep seventh place unchanged. A bob alters this so that the thirteenth change keeps the third place unchanged instead. A single alters this so that the thirteenth change keeps the third place unchanged, and the fourteenth change keeps the first, second and third places all unchanged.

In the case of both Stedman and Erin Triples, for every sixth change, the seventh place remains unchanged. In both cases, a bob alters this so that the sixth change instead keeps the fifth place unchanged. A single makes it so that in addition to the seventh place remaining unchanged in the sixth change, the fifth and sixth places also remain unchanged. In the case of Stedman triples, it is necessary to designate whether the call occurs for the sixth change, or the twelfth change. To handle this, the twelve changes are divided into two sets of six changes, with the first six called {\em slow}, and the second six called {\em quick}. In the case where no alternative is chosen, the sixth change is said to be a {\em plain}. Note that a plain is not considered to be a call.

A fundamental question in change ringing is whether there exists peals that follow certain methods and only use particular calls. The question of whether it was possible to use only bobs and plains to obtain a peal for triples was first considered in the 1700s. Since plains are not considered calls, such peals are called {\em bobs-only}. The earliest result along these lines was due to Thompson \cite{thompson} in 1886 who proved that it is impossible to construct a full peal of bobs-only Grandsire Triples. This result was later generalised by Rankin \cite{rankin,rankin2}, and Swan \cite{swan} subsequently provided a simplified proof of Rankin's result. In the late 1980s, White discussed Stedman and Erin Triples at length in a series of two papers \cite{white,white2}. The problem for Stedman Triples remained open until 1994 when Wyld constructed a peal of bob-only Stedman Triples, but did not publish it at the time. Johnson and Saddleton \cite{saddleton} subsequently found and published a solution online in 1995, after which Wyld's peal was also published \cite{wyld}. Johnson \cite{johnson} later found an improved peal (in the sense that it has a simpler and sparser pattern of calls to remember), and in 2012 improved on it again. The problem remains open for Erin Triples despite significant attempts to resolve it.

In this manuscript, we investigate bobs-only Stedman and Erin Triples in the framework of a famous problem from graph theory, the Hamiltonian cycle problem (HCP). Specifically, we show that the search for bobs-only Triples can be thought of as a constrained form of HCP. We then take advantage of a recently discovered family of subgraphs that allows us to convert this problem to a standard instance of HCP, so that this classical problem can be tackled using the highly sophisticated algorithms available for HCP. We show that the graphs so produced are exceptionally difficult to solve despite their fairly modest size. We also investigate the use of a group theoretic approach to decompose the problem into smaller parts and show that the HCP conversion remains valid.

\section*{Hamiltonian Cycle Problem}

The Hamiltonian cycle problem (HCP) is a, now classical, decision problem from graph theory which can be stated succinctly: given a graph $G$ containing vertex set $V$ with cardinality $N$, and a (directed) edge set $E : V \rightarrow V$, does the graph contain any simple cycles of length $N$? Such simple cycles are called {\em Hamiltonian cycles} (HC), and any graph containing at least one HC is said to be {\em Hamiltonian}, while those graphs with no HCs are called {\em non-Hamiltonian}. HCP has gained notoriety recently because it was one of the first problems identified to be NP-complete \cite{karp}, and because of its close relationship with the arguably more famous traveling salesman problem (TSP). Indeed, TSP can be thought of as the problem of finding the HC (if one exists) of optimal length in a weighted graph.

Although HCP is an NP-complete problem, there is a wealth of sophisticated heuristics available for HCP, most notably LKH \cite{helsgaun} and Snakes and Ladders Heuristic \cite{slh}. The Concorde TSP solver \cite{concorde} is also very effective on HCP instances despite being primarily designed for TSP. There is also an algorithm by Chalaturnyk \cite{chalaturnyk} that enumerates {\em every} HC in a graph, while Eppstein \cite{eppstein} has an enumerative algorithm designed solely for 3-regular graphs, that is, graphs where each vertex is connected to exactly three other vertices. Note that although HCP is defined for directed graphs, all of the best current heuristics are designed for undirected graphs. If necessary, a directed graph can be converted to an undirected graph but the process involves tripling the number of vertices (e.g. see page 46 in Ejov et al. \cite{3hcp}).

Given the definition of HCP, it is natural to think of the search for a peal in the context of HCP since, in both cases, we are seeking a path (method) that visits every vertex (row) precisely once and then returns to the first choice again. In particular, imagine a graph with $7! = 5040$ vertices, where each vertex corresponds to a row, that is, a permutation of seven numbers. Then, for any valid change that can be made for this row, which results in a new row, an edge can be introduced connecting this vertex to the vertex corresponding to the new row. What results is a graph in which the full set of peals and the full set of Hamiltonian cycles are in 1-to-1 correspondence.

If we are searching for a peal that follows a particular method, for example a peal of bobs-only Stedman Triples, then not just any HC will be sufficient. Instead, we need a HC that satisfies the requirements of the method. In the case of bobs-only peals, it is sufficient to only consider changes where a single place is unchanged in the row. Recall that there are only four such changes possible for any given row, so the graph can be reduced to a 4-regular graph. However, even in this reduced graph, the number of possible HCs is likely to be astronomically higher than those corresponding to Stedman Triples. Instead, what we are looking for is a {\em constrained Hamiltonian cycle}, one that follows the pattern of 3.1.3.1.3.(7/5).1.3.1.3.1.(7/5), where the (7/5) changes indicate that either 7th place remains unchanged (a plain) or 5th place remains unchanged (a bob). No standard HCP heuristic is designed to handle such constraints, so the only remaining alternative is to modify the graph in such a way that HCs not satisfying this constraint are eliminated, while those that do are retained.

Before we do this, however, it is possible to reduce the problem further. Since the only decisions that can be made occur every six steps, it seems sensible to consider groups of six changes together. In the case of Stedman and Erin Triples, this is greatly facilitated by the existence of {\em sixes}, which are sets of six rows. The full extent of 5040 rows can be partitioned uniquely into 840 sixes using the following procedure. First, note that following the pattern 3.1.3.1.3.1 causes the row to repeat. For example, starting from rounds, we go through: 1234567, 2135476, 2314567, 3215476, 3124567, 1325476, and then back to rounds again. Hence, these six rows are considered to be members of the same {\em six}. Irrespective of which of them we start with first, the same six members will be produced. Then we can select any row not yet considered and generate its six, and so on until all 840 sixes are produced. Now, consider the case where a plain or a bob occurs, resulting in a certain row. Then it is clear that the subsequent rows will also be a member of that row's six until the next plain or a bob occurs. Hence, for a bobs-only peal of Stedman Triples or Erin Triples, we can think of the problem as the task of visiting all 840 sixes without repetition and returning to the first six; this version of the problem also naturally lends itself to HCP.

With the consideration of sixes, some new terminology is required. Although each six may only be visited once, the first row visited in any given six is not fixed, and will change depending on the calls made to that stage. Hence, the term {\em six-end} is used to denote the final row visited in a six for a particular peal. In the case of Erin Triples, the six-end is uniquely determined by which row is visited first in the six; upon visiting that row, the pattern 3.1.3.1.3 is traversed resulting in the six-end. In the case of Stedman Triples, upon entering a six, the pattern may be 3.1.3.1.3 (slow) or 1.3.1.3.1 (quick). To represent this, we say that a six can be either a {\em slow six} or a {\em quick six}, the determination of which depends on the peal. Each quick six is followed by a slow six, and vice versa. If a six is concluded by a plain, we say it is a {\em plain six}, and if concluded by a bob it is called a {\em bobbed six}, with the determination again depending on the peal.

Fortunately, a parity argument can be utilised to reduce the size of the problem by half. Each row can be labelled as {\em odd} or {\em even}, such that even rows only link to odd rows and vice versa. Hence, if we start with an even row, the six-end will be an odd row, and hence the new six (whichever it is) will be entered at an even row. Inductively, we can conclude that only odd rows may be six-ends. In the case of Erin Triples, this means that each six has three possible six-ends. In the case of Stedman Triples, each six-end may be slow or quick, so for simplicity we think of them as six separate six-ends.

Subsequently, we can now cast the problem as a smaller version of constrained HCP. First consider Erin Triples. Then we have a directed graph with 840 vertices, each corresponding to a six. From each vertex, there are six outgoing edges, corresponding to the two possible ways (plain or bob) to leave from each of the three six-ends. There will also be six incoming edges from other sixes. Then any peal of bobs-only Erin Triples will correspond to a HC in this graph that satisfies the constraints that, upon entering a vertex (six), one of the two \lq\lq valid" outgoing edges (corresponding to the appropriate six-end) are used. Note that determining which six-end is the appropriate one depends on the previous six-end. Hence, the constraints to be satisfied can all be stated in the form, \lq\lq upon visiting this vertex via edge $e_1$, you may only depart via edges $e_2$ or $e_3$."

In the case of Stedman Triples we also have a directed graph with 840 vertices, each corresponding to a six. However, this time from each vertex there are twelve outgoing edges, since in this case we also need to consider whether the the six is slow or quick. Nonetheless, the nature of the constraints is still the same; upon entering a vertex, there will only be two \lq\lq valid" outgoing edges, and which two are valid depends entirely on the edge used to enter.

Since no HCP algorithm is designed to handle such constraints, we take advantage of subgraphs which we called {\em in-out subgraphs} which will be used to eliminate any HCs that don't satisfy the desired constraints.

\section*{In-out Subgraphs}

A recent paper \cite{brokencrown} introduced the family of {\em Crown subgraphs} $\mathcal{C}_n$, which satisfy the following property, which we call the \lq\lq in-out" property. Suppose that $\mathcal{C}_n$ is included in any large graph $G$ such that:

\begin{enumerate}\item[(1)] There are $n$ vertices in $\mathcal{C}_n$ denoted {\em incoming vertices}, such that any incoming edges from the rest of $G$ to $\mathcal{C}_n$ must go to one of these vertices. These are labelled $i_1, i_2, \hdots, i_n$.
\item[(2)] Similarly, there are $n$ vertices in $\mathcal{C}_n$ denoted {\em outgoing vertices}, such that any outgoing edges from $\mathcal{C}_n$ to the rest of $G$ must depart from one of these vertices. Note that these may overlap with the incoming vertices. These are labelled $o_1, o_2, \hdots, o_n$.\end{enumerate}

Then, consider any HC $H$ contained in $G$. The following must be true:

\begin{enumerate}\item[(3)] Upon entering $\mathcal{C}_n$, the Hamiltonian cycle $H$ must traverse the entirety of $\mathcal{C}_n$ before departing. That is, $\mathcal{C}_n$ cannot be entered and exited twice in any valid HC, irrespective of the structure of the rest of the graph.
\item[(4)] If $H$ enters $\mathcal{C}_n$ via incoming vertex $i_k$, it must depart via outgoing vertex $o_k$.\end{enumerate}

Property (3) ensures that $\mathcal{C}_n$ behaves like a vertex in the context of HCP, in the sense that each copy of $\mathcal{C}_n$ included in a graph must be visited precisely once. Property (4) can be used to impose the kinds of constraint mentioned in the previous section. In particular, if every vertex is replaced by $\mathcal{C}_n$, then depending on which edge is used, a particular incoming vertex will be visited, and this will determine the outgoing vertex and hence which outgoing edges may be used. Hence, by linking 840 copies of $\mathcal{C}_n$ together appropriately, only those peals corresponding to the desired method will remain. In the case of Erin Triples, $\mathcal{C}_3$ can be used, while for Stedman Triples, $\mathcal{C}_6$ can be used.

However, $\mathcal{C}_6$ contains 20 vertices, and since the resulting graph is directed, converting it to an undirected graph (so heuristics can consider it) will triple the size of each $\mathcal{C}_6$ to 60 vertices. This would result in a graph containing 50,400 vertices for Stedman Triples. Meanwhile, $\mathcal{C}_n$ is only defined for $n \geq 4$, so $C_3$ does not exist, although a stripped down version of $\mathcal{C}_4$ can be produced which contains 8 vertices, which upon conversion to an undirected graph would result in a graph containing 20,160 vertices for Erin Triples. In order to improve the solving time for these graphs, it makes sense to search for other, smaller, graphs satisfying the in-out property to be used instead of $\mathcal{C}_n$. Such graphs do not appear to have been widely studied or even considered in literature beyond the results of \cite{brokencrown}, although an example of a 2-in-out subgraph can be seen in Papadimitriou and Steiglitz \cite{papadimitriou}. We will refer to such subgraphs as $n$-in-out subgraphs, where $n$ denotes the number of incoming (and hence outgoing) vertices. For the purposes of considering Erin and Stedman Triples, we have successfully identified a 3-in-out subgraph and a 6-in-out subgraph which are both considerably smaller in size than the equivalent Crown subgraphs.

Before we provide the in-out subgraphs, note that although the graphs considered so far have been directed, since we will be replacing {\em every} vertex with an in-out subgraph, it is sufficient to consider strictly undirected in-out subgraphs, as long as we add the caveat that incoming and outgoing vertices may no longer overlap. Then as long as the HC is started at an incoming vertex (any vertex may be chosen as the starting point of a HC), this will be equivalent to the directed graph.

We denote the two (undirected) in-out subgraphs we discovered as $S_3$ and $S_6$, on 16 and 33 vertices respectively. Both subgraphs can be constructed by starting with a path graph (that is, a graph of $n$ vertices containing edges $(i,i+1)$ for $i = 1, \hdots, n-1$) and then adding the following (undirected) edges.

For $S_3$ : Add edges $(1,5)$, $(1,13)$, $(3,12)$, $(4,8)$, $(9,16)$.\\
For $S_6$ : Add edges $(1,24)$, $(4,9)$, $(6,31)$, $(7,12)$, $(10,15)$, $(13,18)$, $(16,27)$, $(22,33)$, $(25,30)$, $(28,33)$.

Then, for $S_3$ the incoming vertices are $i_1 = 1$, $i_2 = 6$ and $i_3 = 8$, while the outgoing vertices are $o_1 = 16$, $o_2 = 11$, $o_3 = 14$. For $S_6$ the incoming vertices are $i_1 = 1$, $i_2 = 7$, $i_3 = 13$, $i_4 = 19$, $i_5 = 25$ and $i_6 = 31$, while the outgoing vertices are $o_1 = 33$, $o_2 = 27$, $o_3 = 3$, $o_4 = 9$, $o_5 = 21$ and $o_6 = 15$. Although we omit the full proof here, it can be checked by exhaustive search that both $S_3$ and $S_6$ satisfy the in-out property, with the following paths from each incoming vertex to the appropriate outgoing vertex.

For $S_3$:

1 - 2 - 3 - 4 - 5 - 6 - 7 - 8 - 9 - 10 - 11 - 12 - 13 - 14 - 15 - 16\\
6 - 7 - 8 - 4 - 5 - 1 - 2 - 3 - 12 - 13 - 14 - 15 - 16 - 9 - 10 - 11\\
8 - 7 - 6 - 5 - 4 - 3 - 2 - 1 - 13 - 12 - 11 - 10 - 9 - 16 - 15 - 14

For $S_6$:

{\small 1-2-3-4-5-6-7-8-9-10-11-12-13-14-15-16-17-18-19-20-21-22-23-24-25-26-27-28-29-30-31-32-33\\
7-8-9-10-11-12-13-14-15-16-17-18-19-20-21-22-23-24-1-2-3-4-5-6-31-32-33-28-29-30-25-26-27\\
13-14-15-10-11-12-7-8-9-4-5-6-31-32-33-28-29-30-25-26-27-16-17-18-19-20-21-22-23-24-1-2-3\\
19-20-21-22-23-24-1-2-3-4-5-6-31-32-33-28-29-30-25-26-27-16-17-18-13-14-15-10-11-12-7-8-9\\
25-26-27-28-29-30-31-32-33-22-23-24-1-2-3-4-5-6-7-8-9-10-11-12-13-14-15-16-17-18-19-20-21\\
31-32-33-28-29-30-25-26-27-16-17-18-19-20-21-22-23-24-1-2-3-4-5-6-7-8-9-10-11-12-13-14-15}

\begin{figure}[h!]\begin{center}\includegraphics[scale=0.25]{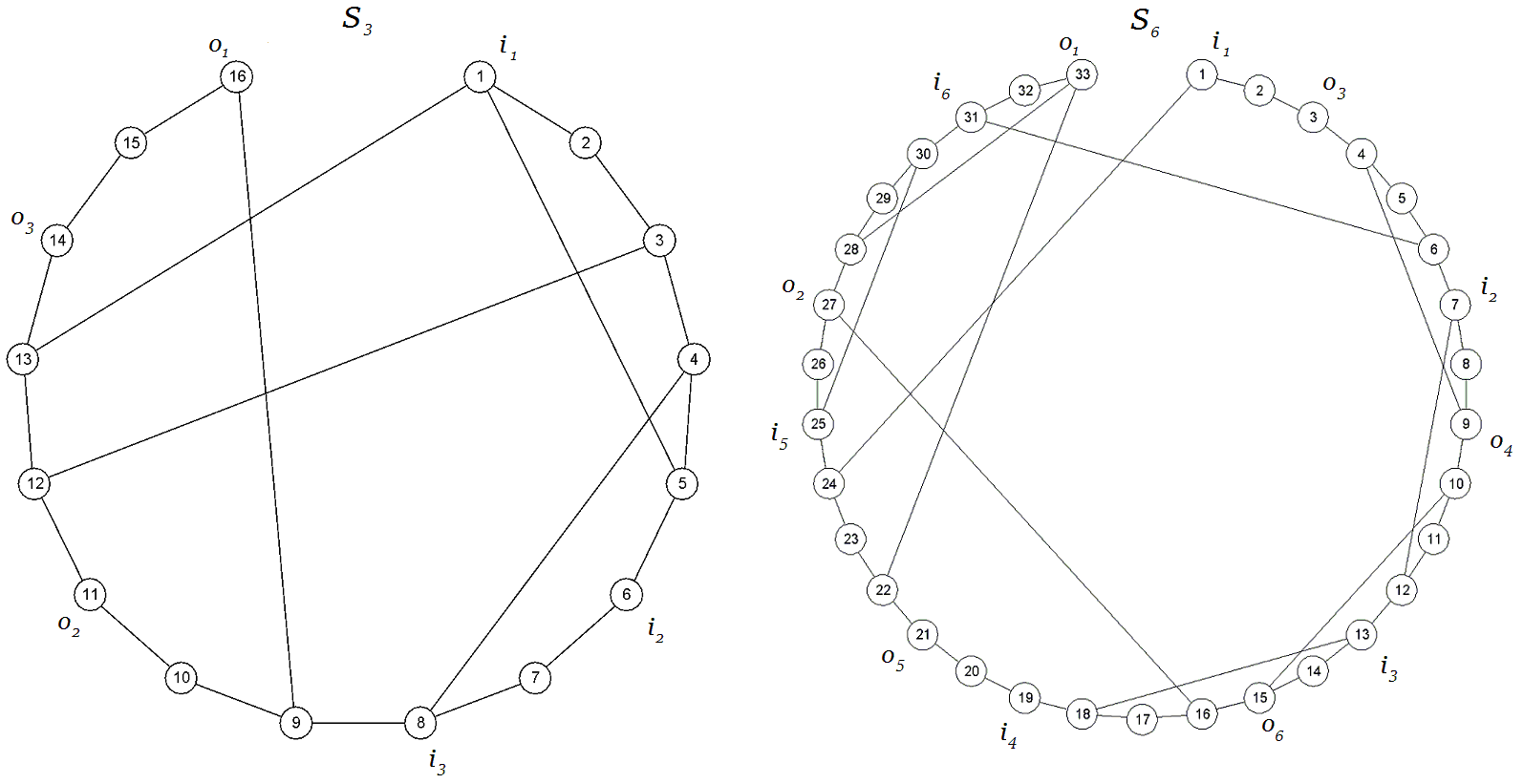}\caption{The in-out subgraphs $S_3$ and $S_6$.}\label{fig-ios}\end{center}\end{figure}

Graphs $S_3$ and $S_6$ are displayed in Figure \ref{fig-ios}.

Using copies of $S_3$, we can produce an instance of HCP corresponding to bobs-only Erin Triples via the following algorithm.

\begin{enumerate}\item Find all 840 sixes, and identify the three six-ends in each six. Denote the sixes $s_1, \hdots, s_{840}$, and each of the three six ends in six $j$ by $e^1_j, e^2_j, e^3_j$.
\item For each six-end $e^k_j$, work out which six-end will be visited next after a plain (ie after going through changes 7.3.1.3.1.3), and which six-end will be visited next after a bob (ie after going through changes 5.3.1.3.1.3). Denote a map $P : e \rightarrow e$ such that $P(e^k_j)$ is equal to the next six-end visited after a plain, and an equivalent map $B$ for a bob.
\item Build a graph from 840 copies of $S_3$. Denote by $S_3^j$ the $j$-th copy of $S_3$.
\item For all $j = 1, \hdots, 840$ and $k = 1, 2, 3$, consider $e_j^k$ and $e_a^b = P(e_j^k)$. Add an edge from outgoing vertex $o_k$ of $S_3^j$ to incoming vertex $i_b$ of $S_3^a$.
\item For all $j = 1, \hdots, 840$ and $k = 1, 2, 3$, consider $e_j^k$ and $e_a^b = B(e_j^k)$. Add an edge from outgoing vertex $o_k$ of $S_3^j$ to incoming vertex $i_b$ of $S_3^a$.\end{enumerate}

Similarly, using copies of $S_6$, we can produce an instance of HCP corresponding to bobs-only Stedman Triples via the following algorithm.

\begin{enumerate}\item Find all 840 sixes, and identify the three six-ends in each six. Denote the sixes $s_1, \hdots, s_{840}$, and each of the three six ends in six $j$ by $e^1_j, e^2_j, e^3_j$.
\item For each six-end $e^k_j$, work out which six-end will be visited next after a slow plain (7.1.3.1.3.1), a quick plain (7.3.1.3.1.3), a slow bob (5.1.3.1.3.1) and a quick bob (5.3.1.3.1.3). Denote the respective maps by $SP$, $QP$, $SB$ and $QB$.
\item Build a graph from 840 copies of $S_6$. Denote by $S_6^j$ the $j$-th copy of $S_6$.
\item For all $j = 1, \hdots, 840$ and $k = 1, 2, 3$, consider $e_j^k$ and $e_a^b = SP(e_j^k)$. Add an edge from outgoing vertex $o_k$ of $S_6^j$ to incoming vertex $i_{b+3}$ of $S_6^a$.
\item For all $j = 1, \hdots, 840$ and $k = 1, 2, 3$, consider $e_j^k$ and $e_a^b = QP(e_j^k)$. Add an edge from outgoing vertex $o_{k+3}$ of $S_6^j$ to incoming vertex $i_b$ of $S_6^a$.
\item For all $j = 1, \hdots, 840$ and $k = 1, 2, 3$, consider $e_j^k$ and $e_a^b = SB(e_j^k)$. Add an edge from outgoing vertex $o_k$ of $S_6^j$ to incoming vertex $i_{b+3}$ of $S_6^a$.
\item For all $j = 1, \hdots, 840$ and $k = 1, 2, 3$, consider $e_j^k$ and $e_a^b = QB(e_j^k)$. Add an edge from outgoing vertex $o_{k+3}$ of $S_6^j$ to incoming vertex $i_b$ of $S_6^a$.\end{enumerate}

The resultant Erin graph has 13,440 vertices and 21,840 edges, while the Stedman graph has 27,720 vertices and 45,360 edges. In the case of the Erin graph, any HC (or a proof that the graph is non-Hamiltonian) would settle the centuries-old question of whether a peal of bobs-only Erin Triples exists.

\section*{Peals in Parts}

One of the breakthroughs in the search for peals was obtained by utilising techniques from group theory to break the full extent up into several equivalent parts. Then, any method that fully traverses one of the parts can be used for the other parts as well, in the sense that the starting six from one part corresponds to a particular six in each of the other parts respectively, and following the same calls, each of the parts can be individually traversed. Each part has a designated starting point, and the final step of the path in each part either returns to this starting point, completing a cycle, or goes to the starting point of another part. Since the parts of equivalent, the latter case means that the path through this new part also finishes by going to the starting point of another part (possibly the first one) and so the final result of each of these paths is a set of short cycles, the union of which forms a vertex-disjoint cycle cover in the graph. These short cycles which contain one or more parts are called {\em round blocks}. Then, by the use of some clever switching techniques (e.g. see Price \cite{price2}), it is possible to replace some calls in a systematic way to start joining the round blocks together. These techniques always increase or decrease the number of round blocks by an even amount (or keep them the same) so the search for a single peal in this manner must start from an odd number of round blocks, and so usually from an odd number of parts; however, searching for peals in an even number of parts is of interest in its own right despite typically being unable to form a full peal from the results.

In practice, the process of searching for peals in parts is done by reducing the set of 840 sixes to some fraction of the sixes by mapping each of the omitted sixes to one of the remaining sixes. Likewise, each of the six-ends in the omitted sixes maps to the equivalent six-end of the corresponding remaining six. What results is a smaller problem which can be easily converted to an instance of HCP using the algorithms in the previous section.

In order to perform this mapping, an appropriate group is required. These cannot be constructed at will, but various groups are known which can be used for either Stedman or Erin triples. The group contains a set of six-ends which will be contained in each of the respective parts. For example, the group used to search for a 5-part peal of bobs-only Stedman Triples is: $\{1234567, 2345167, 3451267, 4512367, 5123467\}$. This group can be produced from the generator $\{2345167\}$. From there, the process to follow is to identify the other rows for each of the five sixes by going through the standard 3.1.3.1.3 process:

$\begin{array}{|c|c|c|c|c|}\hline \;\;\; 1234567 \;\;\; & \;\;\; 2345167 \;\;\; & \;\;\; 3451267 \;\;\; & \;\;\; 4512367 \;\;\; & \;\;\; 5123467 \;\;\;\\
2135476 & 3241576 & 4352176 & 5413276 & 1524376\\
2314567 & 3425167 & 4531267 & 5142367 & 1253467\\
3215476 & 4321576 & 5432176 & 1543276 & 2154376\\
3124567 & 4235167 & 5341267 & 1452367 & 2513467\\
1325476 & 2431576 & 3542176 & 4153276 & 5214376\\
\hline\end{array}$

Note that, although we can move down the column by going through the process 3.1.3.1.3, we can almost move across rows using a relabelling method. For example, from the first entry in the second column to the first entry in the second column, we relabel $\{1,2,3,4,5,6,7\} \rightarrow \{2,3,4,5,1,6,7\}$ so 1 is relabelled to 2, 2 is relabelled to 3, and so on. This relationship is retained throughout the table. Hence, when we next choose any even row that has not been already considered (for example, by considering a plain or bob on one of the six-ends obtained so far) and add it to the first column to obtain a new six, we can use this same relabelling technique to obtain the rest of the row, and proceed to find the next sixes. For example, in the above table, 5324716 had not yet been considered, so we can continue:

$\begin{array}{|c|c|c|c|c|}\hline \;\;\; 1234567 \;\;\; & \;\;\; 2345167 \;\;\; & \;\;\; 3451267 \;\;\; & \;\;\; 4512367 \;\;\; & \;\;\; 5123467 \;\;\;\\
2135476 & 3241576 & 4352176 & 5413276 & 1524376\\
2314567 & 3425167 & 4531267 & 5142367 & 1253467\\
3215476 & 4321576 & 5432176 & 1543276 & 2154376\\
3124567 & 4235167 & 5341267 & 1452367 & 2513467\\
1325476 & 2431576 & 3542176 & 4153276 & 5214376\\
\hline
5324716 & 1435726 & 2541736 & 3152746 & 4213756 \\
3527461 & 4137562 & 5241763 & 1357264 & 2417365 \\
3254716 & 4315726 & 5427136 & 1532746 & 2143756 \\
2357461 & 3417562 & 4521763 & 5137264 & 1247365 \\
2534716 & 3145726 & 4257136 & 5312746 & 1423756 \\
5237461 & 1347562 & 2451763 & 3517264 & 4127365 \\
\hline
\vdots & \vdots & \vdots & \vdots & \vdots\\
\hline
\end{array}$

It can be shown (e.g. see Price \cite{price}) that it is impossible to obtain any repetition this way, and hence this procedure results in a partitioning of the 840 sixes into five disjoint parts of 168 sixes each. Then, each six in columns 2, 3, 4 or 5 are just mapped to equivalent the six in column 1. Likewise, any six-end in columns 2, 3, 4 or 5 is mapped to the equivalent six-end in column 1. So, for example, consider the six-end 1325476 in the first column. After a plain, the next row would be 3152746. Where the next six-end occurs depends on whether 1325476 was a slow six-end or a quick six-end. To represent the different cases, we write either S or Q at the end. Therefore, the six-end after choosing a plain from 1325476S is 1357264Q, and the six-end after choosing a plain from 1325476Q is 3517264S. Both of these entries are in the second six in the fourth column, so the would be mapped to the first column accordingly. That is, a plain from 1325476S would be mapped to 3527461Q, and a plain from 1325476Q would be mapped to 5237461S.

Once the entire set of mappings has been created, the algorithm in the previous section can be used to produce the 5-part Stedman graph. Any time a six-end is mapped to another six-end in the same six, the edge can be ignored, making the graph a bit simpler (although no issue arises if the edge is retained). The resulting graph is one-fifth the size of the original Stedman graph, making it more tractible for HCP heuristics.

Price \cite{price2} gave a large set of groups that have been discovered which are useful for change ringing problems. Here we summarise those that could, potentially, be useful for finding bobs-only Stedman or Erin Triples. They are selected by first considering the order of the group (a $k$-part peal requires a group of order $k$) and then choosing only those groups which satisfy the following conditions. First, there must not be a 3-bell cycle within one of the elements. Second, there must not be 3 pairs swapping in one of the elements. These conditions ensure the splitting into parts does not interfere with the splitting into sixes. For bobs-only peals we also want an odd number of round blocks, so the group must have a cyclic subgroup with an odd number of cosets. This cyclic subgroup must also be even so the start of each part is reachable with just bobs and plains. In each of the remaining cases we give the generators, as well as the index used by Price \cite{price}, and the potential numbers of round blocks that might occur.

1-part: Index 0.01: $\{1234567\}$. 1 round block.\\
2-part: Index 4.07: $\{2143567\}$. 1 or 2 round blocks.\\
3-part: Index 6.33: $\{1357246\}$. 1 or 3 round blocks.\\
4-part: Index 6.26: $\{1352476\}$. 1, 2 or 4 round blocks.\\
5-part: Index 5.05: $\{2345167\}$. 1 or 5 round blocks.\\
6-part: Index 6.32: $\{2316457, 2136547\}$. 2, 3 or 6 round blocks.\\
7-part: Index 7.07: $\{2345671\}$. 1 or 7 round blocks.\\
10-part: Index 5.04: $\{2345167, 5432167\}$. 2, 5 or 10 round blocks.\\
20-part: Index 7.12: $\{2345167, 1352476\}$. 4, 5, 10 or 20 round blocks.\\
21-part: Index 7.05: $\{2345671, 1357246\}$. 3, 7 or 21 round blocks.

\section*{Results}

For each of the ten groups mentioned in the previous section, we produced an HCP instance for both Stedman and Erin Triples. Using exhaustive search techniques, we were able to either discover solutions for each of the Stedman cases, or confirm that none existed in the case of Sted21 and Sted7. In many cases, the solution for a smaller graph, when translated back to a peal, could then be utilised to provide solutions for larger graphs, for example a peal arising from Sted10 can also be seen as arising from an equivalent solution of Sted5 in the sense that the latter is just a repeated version of the former. For the Erin cases, we were able to confirm that none exists in all cases except for the 2-part and 1-part graphs. For the 2-part graph, we found Hamiltonian cycles, but none could be used to find a valid peal. We have not yet been able to confirm if a solution exists for the 1-part graph, but the existence of a solution for the 2-part graph provides some hope. For the Stedman graphs up to and including Sted5, we were able to determine the precise number of solutions. It is worth noting that the conversion from Stedman and Erin Triples to HCP is a one-to-one conversion in the sense that there is precisely one Hamiltonian cycle corresponding to every valid peal. These results are summarised in Table \ref{tab-results}.

\begin{table}[h!]\begin{center}\begin{tabular}{|l|l|l|l|l|l|l|l|l|}
\hline {\bf Graph} & {\bf Nodes} & {\bf Edges} & {\bf Hamiltonicity} & & {\bf Graph} & {\bf Nodes} & {\bf Edges} & {\bf Hamiltonicity}\\
\hline Sted1 & 27720 & 45360 & H & & Erin1 & 13440 & 21840 & U\\
\hline Sted2 & 13860 & 22680 & H & & Erin2 & 6720 & 10920 & H\\
\hline Sted3 & 9240 & 15120 & H  & & Erin3 & 4480 & 7280 & NH\\
\hline Sted4 & 6930 & 11340 & H & & Erin4 & 3360 & 5460 & NH\\
\hline Sted5 & 5544 & 9072 & H (4)& & Erin5 & 2688 & 4368 & NH\\
\hline Sted6 & 4620 & 7560 & H (132) & & Erin6 & 2240 & 3640 & NH\\
\hline Sted7 & 3960 & 6480 & NH & & Erin7 & 1920 & 3120 & NH\\
\hline Sted10 & 2772 & 4536 & H (4) & & Erin10 & 1344 & 2184 & NH\\
\hline Sted20 & 1386 & 2268 & H (6) & & Erin20 & 672 & 1092 & NH\\
\hline Sted21 & 1320 & 2160 & NH & & Erin21 & 640 & 1040 & NH\\
\hline\end{tabular}\caption{The 20 graphs produced for Stedman and Erin Triples, as well as their Hamiltonicity, where 'U' indicates the Hamiltonicity is currently unknown. For some Hamiltonian graphs, the number of solutions is given in brackets after the 'H' marking.}\label{tab-results}\end{center}\end{table}

The Stedman graphs which are known to be Hamiltonian, as well as Erin2, provide exceptionally difficult instances of HCP. This difficulty is interesting in its own right, since truly difficult instances of HCP are not widely known in literature and there has not been a truly difficult benchmark set produced to date. Recently, an HCP challenge was held, called the FHCP Challenge \cite{fhcpcs}, to take strides in this direction. In the FHCP Challenge, which ran for one year, a set of 1001 graphs known to be Hamiltonian needed to be solved. The first team to solve the graphs would win the prize of \$1001US, and if no teams were successful after a year, the team with the most correctly solved graphs would win. At the conclusion of the FHCP Challenge, only 16 of the graphs remained unsolved by any of the submissions, all of which were based on the (Hamiltonian) Stedman graphs. Specifically, the graphs Sted4, Sted5, Sted6, as well as some graphs that combined them with other graphs\footnote{Sted1 and Sted2 were not included in the FHCP Challenge as they were considered too large for the challenge, while the Hamiltonicity of Erin2 and Sted3 was not yet determined at the time of the challenge.}. This result highlights not only the enormous difficulty of the underlying problem of finding Stedman Triples, but also the lack of other truly difficult instances of HCP of moderate size. It is possible that other difficult instances of HCP might be produced by considering peals in parts for $n \geq 8$ bells by extending the techniques used in this manuscript, which could be extremely useful for stress-testing future HCP heuristics.

Given the difficulty of finding solutions even to Sted10, which takes several days of computational time on standard heuristics such as Snakes and Ladders Heuristic \cite{slh} or Concorde \cite{concorde}, there is still hope that a solution exists for Erin1 despite one not having been found so far. If so, the continued evolution of HCP heuristics and parallel computing technology should make finding a HC in Erin1, and hence a bobs-only peal of Erin Triples, a more tractable problem.

Finally, if the desire is to create difficult instances of HCP, rather than find bobs-only peals, groups which do not satisfy the condition of separating into an odd number of round blocks can be used. These are as follows:

4-part: Index 4.04: $\{2143567, 3412567\}$. 2 or 4 round blocks.\\
4-part: Index 6.35: $\{2143567, 2134657\}$. 2 or 4 round blocks.\\
8-part: Index 6.23: $\{2341657, 4321567\}$. 2, 4 or 8 round blocks.\\
12-part: Index 6.14: $\{5146237, 6423157\}$. 4, 6 or 12 round blocks.\\
12-part: Index 7.33: $\{3476521, 7514623\}$. 4, 6 or 12 round blocks.\\
24-part: Index 6.09: $\{3456127, 2165347\}$. 6, 8, 12 or 24 round blocks.\\
24-part: Index 7.28: $\{2341657, 2134576\}$. 6, 8, 12 or 24 round blocks.\\
60-part: Index 6.05: $\{2345167, 5342617\}$. 12, 20, 30 or 60 round blocks.\\
168-part: Index 7.03: $\{7613524, 4725163\}$. 24, 42, 56, 84 or 168 round blocks.\\

The Hamiltonicity of many of the resulting graphs is known, with only the Stedman graph corresponding to the 4-part group 4.04 not yet determined. These results are summarised in Table \ref{tab-results2}. Note that the cases of Sted8, Sted24 and Sted168 are trivially non-Hamiltonian as the length of each part is not a multiple of 12, so it is impossible to maintain the alternation of quick and slow sixes.

\begin{table}[h!]\begin{center}\begin{tabular}{|l|l|l|l|l|l|l|l|l|}
\hline {\bf Graph} & {\bf Nodes} & {\bf Edges} & {\bf Hamiltonicity} & & {\bf Graph} & {\bf Nodes} & {\bf Edges} & {\bf Hamiltonicity}\\
\hline Sted4 (4.04) & 6930 & 11340 & U & & Erin4 (4.04) & 3360 & 5460 & NH\\
\hline Sted4 (6.35) & 6930 & 11340 & H & & Erin4 (6.35) & 3360 & 5460 & NH\\
\hline Sted8 (6.23) & 3465 & 5670 & NH  & & Erin8 (6.23) & 1680 & 2730 & NH\\
\hline Sted12 (6.14) & 2310 & 3780 & H (248) & & Erin12 (6.14) & 1120 & 1820 & NH\\
\hline Sted12 (7.33) & 2310 & 3780 & NH & & Erin12 (7.33) & 1120 & 1820 & NH\\
\hline Sted24 (6.09) & 1155 & 1890 & NH & & Erin24 (6.09) & 560 & 910 & NH\\
\hline Sted24 (7.28) & 1155 & 1890 & NH & & Erin24 (7.28) & 560 & 910 & NH\\
\hline Sted60 (6.05) & 462 & 756 & H (20) & & Erin60 (6.05) & 224 & 364 & NH\\
\hline Sted168 (7.03) & 165 & 270 & NH & & Erin168 (7.03) & 80 & 130 & NH\\
\hline\end{tabular}\caption{The 18 additional graphs produced for Stedman and Erin Triples for groups that do not permit an odd number of round blocks, as well as their Hamiltonicity, where 'U' indicates the Hamiltonicity is currently unknown. The number of solutions for Sted12 (6.14) and Sted60 (7.03) is also given.}\label{tab-results2}\end{center}\end{table}

\end{document}